\def\author{Franco, Kleiman, and Lascu}
\def\title{Gherardelli linkage and complete intersections}
\def\date{March 9, 2000}
\def\abstract{
   Our main theorem characterizes the complete intersections of
codimension 2 in a projective space of dimension 3 or more over an
algebraically closed field of characteristic 0 as the subcanonical and
self-linked subschemes.  In order to prove this theorem, we'll prove the
Gherardelli linkage theorem, which asserts that a partial intersection
of two hypersurfaces is subcanonical if and only if its residual
intersection is, scheme-theoretically, the intersection of the two
hypersurfaces with a third.
 }
 \def\PaperSize{letter}        
 \let\@=@  

\def\GetNext#1 {\def\NextOne{#1}\if\relax\NextOne\let\next=\relax
        \else\let\next=\DoIt \fi \next}
\def\DoIt{\Act\NextOne\GetNext}
\def\ActOn#1{\expandafter\GetNext #1\relax\ }
\def\defcs#1{\expandafter\xdef\csname#1\endcsname}

\parskip=0pt plus 1.75pt \parindent10pt
\hsize32pc 
\vsize48pc 
\abovedisplayskip 4pt plus3pt minus1pt
\belowdisplayskip=\abovedisplayskip
\abovedisplayshortskip 2.5pt plus2pt minus1pt
\belowdisplayshortskip=\abovedisplayskip

\def\TRUE{TRUE} 
\def\TheMagstep{\magstep1}
\ifx\DoublepageOutput\TRUE \def\TheMagstep{\magstep0} \fi
\mag=\TheMagstep

\newskip\vadjustskip \vadjustskip=0.5\normalbaselineskip
\def\centertext
 {\hoffset=\pgwidth \advance\hoffset-\hsize
  \advance\hoffset-2truein \divide\hoffset by 2\relax
  \voffset=\pgheight \advance\voffset-\vsize
  \advance\voffset-2truein \divide\voffset by 2\relax
  \advance\voffset\vadjustskip
 }
\newdimen\pgwidth\newdimen\pgheight
\def\letter{letter}\def\AFour{AFour}
\ifx\PaperSize\letter
 \pgwidth=8.5truein \pgheight=11truein
 \message{- Got a paper size of letter.  }\centertext
\fi
\ifx\PaperSize\AFour
 \pgwidth=210truemm \pgheight=297truemm
 \message{- Got a paper size of AFour.  }\centertext
\fi

\def\today{\ifcase\month\or     
 January\or February\or March\or April\or May\or June\or
 July\or August\or September\or October\or November\or December\fi
 \space\number\day, \number\year}
\nopagenumbers
 \newcount\pagenumber \pagenumber=1
 \def\advancepagenumber{\global\advance\pagenumber by 1}
\def\folio{\number\pagenum} 
\headline={%
  \ifnum\pagenum=0\hfill
  \else
   \ifnum\pagenum=1\firstheadline
   \else
     \ifodd\pagenum\oddheadline
     \else\evenheadline\fi
   \fi
  \fi
}
\expandafter\ifx\csname date\endcsname\relax \let\dato=\today
            \else\let\dato=\date\fi
\let\firstheadline\hfill
\def\oddheadline{\eightpoint \rlap{\dato}
 \hfil\headtitle\hfil\llap{\folio}}
\def\evenheadline{\eightpoint\rlap{\folio}
 \hfil\author\hfil\llap{\dato}}
\def\headtitle{\title}

 \newdimen\fullhsize \newbox\leftcolumn
 \def\fulline{\hbox to \fullhsize} 
\def\doublepageoutput
{\let\lr=L
 \output={\if L\lr
           \global\setbox\leftcolumn=\columnbox \global\let\lr=R%
          \else \doubleformat \global\let\lr=L
          \fi
        \ifnum\outputpenalty>-20000 \else\dosupereject\fi
        }%
 \def\doubleformat{\shipout\vbox{%
     \ifx\PaperSize\AFour
           \fulline{\hfil\box\leftcolumn\hfil\columnbox\hfil}%
     \else
           \fulline{\hfil\hfil\box\leftcolumn\hfil\columnbox\hfil\hfil}%
     \fi             }%
     \advancepageno
}
 \def\columnbox{\vbox
   {\if E\topmark\headline={\hfil}\nopagenumbers\fi
    \makeheadline\pagebody\makefootline\advancepagenumber}%
   }%
\fullhsize=\pgheight \hoffset=-1truein
 \voffset=\pgwidth \advance\voffset-\vsize
  \advance\voffset-2truein \divide\voffset by 2
  \advance\voffset\vadjustskip
 
\ifx\FirstPageOnRight\TRUE 
 \null\vfill\nopagenumbers\eject\pagenum=1\relax
\fi
}
\ifx\DoublepageOutput\TRUE \let\pagenum=\pagenumber\doublepageoutput
 \else \let\pagenum=\pageno \fi

 \font\twelvebf=cmbx12          
 \font\smc=cmcsc10              
 \font\tenbi=cmmi14
 \font\sevenbi=cmmi10 \font\fivebi=cmmi7
 \newfam\bifam \def\bi{\fam\bifam} \textfont\bifam=\tenbi
 \scriptfont\bifam=\sevenbi \scriptscriptfont\bifam=\fivebi
 \mathchardef\variablemega="7121 \def\bigomega{{\bi\variablemega}} 
 \mathchardef\variablenu="7117 
\catcode`\@=11          
\def\eightpoint{\eightpointfonts
 \setbox\strutbox\hbox{\vrule height7\p@ depth2\p@ width\z@}%
 \eightpointparameters\eightpointfamilies
 \normalbaselines\rm
 }
\def\eightpointparameters{%
 \normalbaselineskip9\p@
 \abovedisplayskip9\p@ plus2.4\p@ minus6.2\p@
 \belowdisplayskip9\p@ plus2.4\p@ minus6.2\p@
 \abovedisplayshortskip\z@ plus2.4\p@
 \belowdisplayshortskip5.6\p@ plus2.4\p@ minus3.2\p@
 }
\newfam\smcfam
\def\eightpointfonts{%
 \font\eightrm=cmr8 \font\sixrm=cmr6
 \font\eightbf=cmbx8 \font\sixbf=cmbx6
 \font\eightit=cmti8
 \font\eightsmc=cmcsc8
 \font\eighti=cmmi8 \font\sixi=cmmi6
 \font\eightsy=cmsy8 \font\sixsy=cmsy6
 \font\eightsl=cmsl8 \font\eighttt=cmtt8}
\def\eightpointfamilies{%
 \textfont\z@\eightrm \scriptfont\z@\sixrm  \scriptscriptfont\z@\fiverm
 \textfont\@ne\eighti \scriptfont\@ne\sixi  \scriptscriptfont\@ne\fivei
 \textfont\tw@\eightsy \scriptfont\tw@\sixsy \scriptscriptfont\tw@\fivesy
 \textfont\thr@@\tenex \scriptfont\thr@@\tenex\scriptscriptfont\thr@@\tenex
 \textfont\itfam\eightit        \def\it{\fam\itfam\eightit}%
 \textfont\slfam\eightsl        \def\sl{\fam\slfam\eightsl}%
 \textfont\ttfam\eighttt        \def\tt{\fam\ttfam\eighttt}%
 \textfont\smcfam\eightsmc      \def\smc{\fam\smcfam\eightsmc}%
 \textfont\bffam\eightbf \scriptfont\bffam\sixbf
   \scriptscriptfont\bffam\fivebf       \def\bf{\fam\bffam\eightbf}%
 \def\rm{\fam0\eightrm}%
 }
\def\vfootnote#1{\insert\footins\bgroup
 \eightpoint\catcode`\^^M=5\leftskip=0pt\rightskip=\leftskip
 \interlinepenalty\interfootnotelinepenalty
  \splittopskip\ht\strutbox 
  \splitmaxdepth\dp\strutbox \floatingpenalty\@MM
  \leftskip\z@skip \rightskip\z@skip \spaceskip\z@skip \xspaceskip\z@skip
  \textindent{#1}\footstrut\futurelet\next\fo@t}

\def\p.{p.\penalty\@M \thinspace}
\def\pp.{pp.\penalty\@M \thinspace}
\newcount\sctno
\def\sctn#1\par
  {\removelastskip\vskip0pt plus4\normalbaselineskip \penalty-250
  \vskip0pt plus-4\normalbaselineskip \bigskip\medskip
  \centerline{\smc#1}\nobreak\medskip
}

\def\sct#1 {\sctno=#1\relax\sctn#1.\enspace}

\def\item#1 {\par\indent\indent\indent
 \hangindent3\parindent
 \llap{\rm (#1)\enspace}\ignorespaces}
 \def\inpart#1 {{\rm (#1)\enspace}\ignorespaces}
 \def\part {\par\inpart}

\def\Cs#1){\(\number\sctno.#1)}
\def\part#1 {\par\(#1)\enspace\ignorespaces}

\def\dsc#1 #2.{\medbreak{\bf\Cs#1)} {\it #2.} \ignorespaces}
\def\proclaim#1 #2 {\medbreak
  {\bf#1 (\number\sctno.#2)}\enspace \bgroup
\it}
\def\endproclaim{\par\egroup\medskip}
\def\pf{\endproclaim{\bf Proof.} \ignorespaces}
 \def\prp{\proclaim Proposition }
  \def\thm{\proclaim Theorem }
\def\dfn#1 {\medbreak {\bf Definition (\number\sctno.#1)}\enspace}
\def\rmk#1 {\medbreak {\bf Remark (\number\sctno.#1)}\enspace}
\def\stp#1 {\medbreak {\bf Setup (\number\sctno.#1)}\enspace}
\def\eg#1 {\medbreak {\bf Example (\number\sctno.#1)}\enspace}
 \newcount\refno \refno=0        \def\NoKey{*!*}
 \def\MakeKey{\advance\refno by 1 \expandafter\xdef
  \csname\TheKey\endcsname{{\number\refno}}\NextKey}
 \def\NextKey#1 {\def\TheKey{#1}\ifx\TheKey\NoKey\let\next\relax
  \else\let\next\MakeKey \fi \next}
 \def\RefKeys #1\endRefKeys{\expandafter\NextKey #1 *!* }
 \def\SetRef#1 #2,{\hang\llap
  {[\csname#1\endcsname]\enspace}{\smc #2},}
 \newbox\keybox \setbox\keybox=\hbox{[25]\enspace}
 \newdimen\keyindent \keyindent=\wd\keybox
\def\references{\kern-\medskipamount
  \sctn References\par
  \vskip-\medskipamount
  \bgroup   \frenchspacing   \eightpoint
   \parindent=\keyindent  \parskip=\smallskipamount
   \everypar={\SetRef}\par}
\def\endreferences{\egroup}

 \def\serial#1#2{\expandafter\def\csname#1\endcsname ##1 ##2 ##3
        {\unskip\ {\it #2\/} {\bf##1} (##2), ##3.}} 

\def\UThin{\penalty\@M \thinspace\ignorespaces}
\def\(#1){{\let~=\UThin\rm(#1)}}
\def\relaxnext@{\let\next\relax}
\def\cite#1{\relaxnext@
 \def\nextiii@##1,##2\end@{\unskip\space{\rm
        [\SetKey{##1},\let~=\UThin##2]}}%
 \in@,{#1}\ifin@\def\next{\nextiii@#1\end@}\else
 \def\next{{\rm[\SetKey{#1}]}}\fi\next}
\newif\ifin@
\def\in@#1#2{\def\in@@##1#1##2##3\in@@
 {\ifx\in@##2\in@false\else\in@true\fi}%
 \in@@#2#1\in@\in@@}
\def\SetKey#1{{\bf\csname#1\endcsname}}

\catcode`\@=12  

\let\:=\colon \let\ox=\otimes  
  \let\?=\overline

 \def\onto{\to\mathrel{\mkern-15mu}\to}
 
 \def\IP{{\bf P}}  
\def\smashedlongrightarrow{\setbox0=\hbox{$\longrightarrow$}\ht0=1pt\box0}
\def\risom{\buildrel\sim\over{\smashedlongrightarrow}}
\def\smashedlongleftarrow{\setbox0=\hbox{$\longleftarrow$}\ht0=1pt\box0}
\def\lisom{\buildrel\sim\over{\smashedlongleftarrow}}
 \def\lgto{-\mathrel{\mkern-10mu}\to}
 \def\smashedlgto{\setbox0=\hbox{$\scriptstyle\lgto$}\ht0=1.85pt
        \lower1.25pt\box0}
\def\tto{\buildrel\lgto\over{\smashedlgto}}
 \let\vf=\varphi

\def\Act#1{\defcs{c#1}{{\cal#1}}}               
 \ActOn{B C D E F G H I J K L M N O P Q R S }
\def\Act#1{\defcs{#1}{\mathop{\rm#1}\nolimits}} 
 \ActOn{depth }
\def\Act#1{\defcs{c#1}{\mathop{\it#1}\nolimits}}
 \ActOn{Ann Cok Ext Hom Ker Sym }
\def\Act#1{\defcs{#1}{\hbox{\enspace #1\enspace}}} 
 \ActOn{and by for where with on }
\def\Act#1{\defcs{I#1}{{\bf#1}}}                   
 \ActOn{A B D G R S T Z }

\catcode`\@=11
 \def\activeat#1{\csname @#1\endcsname}
 \def\def@#1{\expandafter\def\csname @#1\endcsname}
 {\catcode`\@=\active \gdef@{\activeat}}

\let\ssize\scriptstyle
\newdimen\ex@   \ex@.2326ex

 \def\requalfill{\cleaders\hbox{$\mkern-2mu\mathord=\mkern-2mu$}\hfill
  \mkern-6mu\mathord=$}
 \def\eqfill{$\m@th\mathord=\mkern-6mu\requalfill}
 \def\deffill{\hbox{$:=$}$\m@th\mkern-6mu\requalfill}
 \def\fiberbox{\hbox{$\vcenter{\hrule\hbox{\vrule\kern1ex
     \vbox{\kern1.2ex}\vrule}\hrule}$}}
 \def\Fiberbox{\rlap{\raise0.75pt\fiberbox}\lower0.75pt\fiberbox}

 \font\arrfont=line10
 \def\Swarrow{\vcenter{\hbox{$\swarrow$\kern-.26ex
    \raise1.5ex\hbox{\arrfont\char'000}}}}

 \newdimen\arrwd
 \newdimen\minCDarrwd \minCDarrwd=2.5pc
        
 \def\findarrwd#1#2#3{\arrwd=#3%
  \setbox\z@\hbox{$\ssize\;{#1}\;\;$}%
 \setbox\@ne\hbox{$\ssize\;{#2}\;\;$}%
  \ifdim\wd\z@>\arrwd \arrwd=\wd\z@\fi
  \ifdim\wd\@ne>\arrwd \arrwd=\wd\@ne\fi}
 \newdimen\arrowsp\arrowsp=0.375em      
 \def\findCDarrwd#1#2{\findarrwd{#1}{#2}{\minCDarrwd}
    \advance\arrwd by 2\arrowsp}
 \newdimen\minarrwd 
 \setbox\z@\hbox{$\longrightarrow$} \minarrwd=\wd\z@

 \def\harrow#1#2#3#4{{\minarrwd=#1\minarrwd%
   \findarrwd{#2}{#3}{\minarrwd}\kern\arrowsp
    \mathrel{\mathop{\hbox to\arrwd{#4}}\limits^{#2}_{#3}}\kern\arrowsp}}
 \def@]#1>#2>#3>{\harrow{#1}{#2}{#3}\rightarrowfill}
 \def@>#1>#2>{\harrow1{#1}{#2}\rightarrowfill}
 \def@<#1<#2<{\harrow1{#1}{#2}\leftarrowfill}
 \def@={\harrow1{}{}\eqfill}
 \def@:#1={\harrow1{}{}\deffill}
 \def@ N#1N#2N{\vCDarrow{#1}{#2}\UpDownarrow}
 \def\UpDownarrow{\uparrow\,\Big\downarrow}

\def@'#1'#2'{\harrow1{#1}{#2}\tarrowfill}
 \def\lgTo{\dimen0=\arrwd \advance\dimen0-2\arrowsp
        \hbox to\dimen0{\rightarrowfill}}
 \def\smashedlgTo{\setbox0=\hbox{$\scriptstyle\lgTo$}\ht0=1.85pt
        \lower1.25pt\box0}
 \def\tto{\buildrel\lgTo\over{\smashedlgTo}}
 \def\tarrowfill{\hfil$\tto$\hfil}  

 \def@={\ifodd\row\harrow1{}{}\eqfill
   \else\vCDarrow{}{}\Vert\fi}
 \def@.{\ifodd\row\relax\harrow1{}{}\hfill
   \else\vCDarrow{}{}.\fi}
 \def@|{\vCDarrow{}{}\Vert}
 \def@ V#1V#2V{\vCDarrow{#1}{#2}\downarrow}
\def@ A#1A#2A{\vCDarrow{#1}{#2}\uparrow}
 \def@(#1){\arrwd=\csname col\the\col\endcsname\relax
   \hbox to 0pt{\hbox to \arrwd{\hss$\vcenter{\hbox{$#1$}}$\hss}\hss}}

 \def\squash#1{\setbox\z@=\hbox{$#1$}\finsm@@sh}
\def\finsm@@sh{\ifnum\row>1\ht\z@\z@\fi \dp\z@\z@ \box\z@}

 \newcount\row \newcount\col \newcount\numcol \newcount\arrspan
 \newdimen\vrtxhalfwd  \newbox\tempbox

 \def\innernewdimen{\alloc@1\dimen\dimendef\insc@unt}
 \def\measureinit{\col=1\vrtxhalfwd=0pt\arrspan=1\arrwd=0pt
   \setbox\tempbox=\hbox\bgroup$}
 \def\setinit{\col=1\hbox\bgroup$\ifodd\row
   \kern\csname col1\endcsname
   \kern-\csname row\the\row col1\endcsname\fi}
 \def\findvrtxhalfsum{$\egroup
  \expandafter\innernewdimen\csname row\the\row col\the\col\endcsname
  \global\csname row\the\row col\the\col\endcsname=\vrtxhalfwd
  \vrtxhalfwd=0.5\wd\tempbox
  \global\advance\csname row\the\row col\the\col\endcsname by \vrtxhalfwd
  \advance\arrwd by \csname row\the\row col\the\col\endcsname
  \divide\arrwd by \arrspan
  \loop\ifnum\col>\numcol \numcol=\col%
 \expandafter\innernewdimen \csname col\the\col\endcsname
     \global\csname col\the\col\endcsname=\arrwd
   \else \ifdim\arrwd >\csname col\the\col\endcsname
      \global\csname col\the\col\endcsname=\arrwd\fi\fi
   \advance\arrspan by -1 %
   \ifnum\arrspan>0 \repeat}
 \def\setCDarrow#1#2#3#4{\advance\col by 1 \arrspan=#1
    \arrwd= -\csname row\the\row col\the\col\endcsname\relax
    \loop\advance\arrwd by \csname col\the\col\endcsname
     \ifnum\arrspan>1 \advance\col by 1 \advance\arrspan by -1%
     \repeat
    \squash{\mathop{
     \hbox to\arrwd{\kern\arrowsp#4\kern\arrowsp}}\limits^{#2}_{#3}}}
 \def\measureCDarrow#1#2#3#4{\findvrtxhalfsum\advance\col by 1%
   \arrspan=#1\findCDarrwd{#2}{#3}%
    \setbox\tempbox=\hbox\bgroup$}
 \def\vCDarrow#1#2#3{\kern\csname col\the\col\endcsname
    \hbox to 0pt{\hss$\vcenter{\llap{$\ssize#1$}}%
     \Big#3\vcenter{\rlap{$\ssize#2$}}$\hss}\advance\col by 1}

 \def\setCD{\def\harrow{\setCDarrow}%
  \def\\{$\egroup\advance\row by 1\setinit}
  \m@th\lineskip3\ex@\lineskiplimit3\ex@ \row=1\setinit}
 \def\endsetCD{$\strut\egroup}
 \def\measure{\bgroup
  \def\harrow{\measureCDarrow}%
  \def\\##1\\{\findvrtxhalfsum\advance\row by 2 \measureinit}%
  \row=1\numcol=0\measureinit}
 \def\endmeasure{\findvrtxhalfsum\egroup}

\newbox\CDbox \newdimen\sdim
 \newcount\savedcount   
 \def\CD#1\endCD{\savedcount=\count11%
   \measure#1\endmeasure
   \vcenter{\setCD#1\endsetCD}%
   \global\count11=\savedcount}
 \catcode`\@=\active

 \RefKeys
 AK BP BE B E EC Fa F FL81 FL87 Gh M PS R79 R82
 \endRefKeys
{\leftskip=0pt plus1fill \rightskip=\leftskip
 \obeylines
 \null\bigskip\bigskip\bigskip\bigskip
 {\twelvebf
  \title
 } \bigskip
 \footnote{}{\noindent %
 MSC-class: 14M10 (Primary) 14M06, 1407 (Secondary).}
 \smc Davide Franco
 {\eightpoint\it\medskip
  Dipartimento di Mathematica, Universit\`a di Ferrara
 via Machiavelli {\sl 35}, 44100 Ferrara, Italy
 \rm E-mail: \tt frv\@dns.unife.it \medskip
 }
 Steven L. Kleiman
 {\eightpoint\it\medskip
 Department of Mathematics, Room {\sl 2-278} MIT,
 {\sl77} Mass Ave, Cambridge, MA {\sl02139-4307}, USA
 \rm E-mail: \tt Kleiman\@math.mit.edu \medskip
 } and \medskip
 Alexandru T. Lascu
 {\eightpoint\it\medskip
 Dipartimento di Mathematica, Universit\`a di Ferrara
 via Machiavelli {\sl 35}, 44100 Ferrara, Italy
 \rm E-mail: \tt lsl\@dns.unife.it
 } \medskip
 \rm \dato \bigskip
}
{\advance\leftskip by 1.5\parindent \advance\rightskip by 1.5\parindent
 \eightpoint \noindent
 {\smc Abstract.}\enspace \ignorespaces \abstract \par
}

\sct1 Introduction

Our main result is Theorem (3.2).  It characterizes the complete
intersections of codimension 2 in $\IP^n$ where $n\ge3$, over an
algebraically closed field of characteristic 0, among the
Cohen--Macaulay $X$ as those that are subcanonical and self-linked.
This characterization was formulated by Ellia (pvt.\ comm.), who proved
it in a joint work with Beorchia \cite{BE, Thm.~5, p.~556} assuming $X$
is smooth.  In Remark~6.1 on \p.557, Beorchia and Ellia said they don't
know whether the smoothness ``can be avoided.''  It can!  Furthermore,
$X$ can be reducible and nonreduced.

More precisely, an $X$ is said to be {\it $a$-subcanonical\/} if its
dualizing sheaf $\bigomega_X$ is of the form $\bigomega_X=\cO_X(a)$.  An
$X$ is said to be {\it self-linked\/} by two hypersurfaces $F_1$ and
$F_2$ if $X$ is equal to its own residual scheme in the complete
intersection of $F_1$ and $F_2$.  For example, suppose $X$ is the
complete intersection of $F_1$ and $F_3$.  Then $X$ is self-linked by
$F_1$ and $F_2$ where $F_2:=2F_3$ or where $F_2$ is, more generally, any
hypersurface such that $F_1\bigcap F_2=F_1\bigcap 2F_3$.  Furthermore,
$X$ is $a$-subcanonical where $a$ is the following integer: denote the
degree of $F_i$ by $m_i$; then $a:=m_1+m_3-n-1$.  Now, Theorem (3.2)
says that this is, in fact, the only example!

This second formulation of Theorem (3.2) is more refined than the first.
After all, the first says nothing much about the hypersurfaces $F_i$
involved.  In particular, the first does not suggest anything like the
equation $m_2=2m_3$.  Indeed, in Corollary~4 on \p.557, Beorchia and
Ellia offered an alternative proof of the first formulation in the case
where $X$ is a curve and $m_3\ge m_2\ge m_1$.  The proof is correct, but
the case is vacuous!

Our proof of Theorem (3.2) follows, to a fair extent, the lines of
Beorchia and Ellia's proof of their Theorem 5.  In both proofs, a key
step is to split the normal bundle of $X$ in $\IP^n$.  At this stage, if
$n\ge4$ and $X$ is smooth, then we're done simply because the normal
bundle splits; indeed, Basili and Peskine \cite{BP, p.~87} proved that
then $X$ is a complete intersection.  However, in order to prove Theorem
(3.2) in full generality, we must split the normal bundle with care.
For example, consider the twisted cubic space curve $X$; its normal
bundle is split because $X$ is rational, and it is known that $X$ is
self-linked by a quadric cone and a cubic surface, but $X$ is, of
course, not a complete intersection, nor even subcanonical.

To split the normal bundle, we'll use the Gherardelli linkage theorem
(2.5).  It asserts that, when two hypersurfaces $F_1$ and $F_2$ of
$\IP^n$ intersect partially in an $X$, then $X$ is subcanonical if and
only if its residual scheme $Y$ is, scheme-theoretically, of the form
$Y=F_1\bigcap F_2\bigcap F_3$ where $F_3$ is a suitable hypersurface.
(Such a $Y$ is called a {\it quasi-complete intersection}.)  In
particular, if $X$ is subcanonical, and is self-linked by $F_1$ and
$F_2$, then $X=F_1\bigcap F_2\bigcap F_3$.  In this case, we'll form the
conormal bundles of $X$ in $F_3$ and in $\IP^n$, and we'll split the
natural map from the latter bundle onto the former.

We'll then conclude that some multiple of $X$ is numerically equivalent
to a hypersurface section of $F_3$, at least after we've replaced $F_3$
by an integral component; we'll simply apply Braun's main theorem
\cite{B, p.~403}.  (Braun followed the lines of Ellingsrud, Gruson,
Peskine, and Str\o mme's remarkable proof of the theorem in the case of
a curve on a smooth connected surface. This case had been treated
earlier, in a very different fashion, by Griffiths, Harris, and Hulek.
See Braun's paper \cite{B, p.~411} for all the references.)  Finally, to
conclude that $X$ is a complete intersection, Beorchia and Ellia used
Gruson and Peskine's work on space curves.  Instead, we'll make a direct
geometric argument, and obtain our more refined statement of Theorem
(3.2).

If $n\ge6$ and $X$ is smooth, then, since $X$ is a quasi-complete
intersection, it is, in fact, a complete intersection by Faltings'
Korollar of Satz~3 \cite{Fa, p.~398}.  This line of proof is significant
because it is valid in any characteristic, whereas Basili and Peskine
work in characteristic 0, and we must too although only to apply Braun's
theorem.  Beorchia and Ellia \cite{BE, p.~556} suggested that there
might be a problem in characteristic 2 by pointing out the following
result, due in part to Rao \cite{R82, p.272} and in part to Migliore
\cite{M, p.~185}: a double line in $\IP^3$ of arithmetic genus $-2$ or
less is self-linked if and only if the characteristic is 2.  We'll
pursue this suggestion in Example (3.4).  On the other hand, it would be
nice to know whether Theorem (3.2) is valid except for certain $X$ of
small dimension in characteristic 2.

The Gherardelli linkage theorem holds in greater generality than that
stated above.  In Theorem (2.5), we'll replace $\IP^n$ by any Gorenstein
projective scheme $P$ having pure dimension 2 or more and satisfying
this vanishing condition: $H^q(\cO_P(m))=0$ for three specific values of
the pair $(q,\,m)$.  For example, $P$ can be a complete intersection in
$\IP^n$.  Thus we'll recover Theorem~2(i) of Fiorentini and Lascu
\cite{FL87, p.~170}, where, in addition, $X$ and $Y$ are assumed to have
no common components; in fact, our proof was inspired by theirs.
Beorchia and Ellia \cite{BE, p.~556} proved the existence of the
hypersurface $F_3$ directly in the case at hand by using the mapping
cone.  Earlier, Rao \cite{R79, pp.~209--10} proved (burying it among
other things) a version of the Gherardelli linkage theorem, in which the
condition that $X$ be subcanonical is replaced by the condition that $X$
be the zero scheme of a section of a rank-2 vector bundle on $\IP^n$;
these two conditions are equivalent by a famous theorem of Serre's (see
\cite{F, Prop.~3, p.~346}).  On the other hand, our Theorem (3.2) does
not hold even if $\IP^n$ is replaced by a smooth hypersurface $P$, as
we'll see in Example~(3.3).

To prove the Gherardelli linkage theorem (2.5), we'll use the Noether
linkage sequence (2.3.1), which presents the dualizing sheaf of a
partial intersection in any Gorenstein ambient scheme $P$ having pure
dimension 2 or more.  The case where $P$ is a complete intersection in
$\IP^n$ was treated in \cite{FL87, Lem.~1} and in \cite{PS, 1.6} and was
used in \cite{R82, p.~253}.  The general case is, as we'll see, no more
difficult to prove.

In short, in Section 2, we'll review some basic linkage theory,
including the Peskine--Szpiro linkage theorem (compare with \cite{PS, 1.3,
\p.274} and \cite{E, 21.23, \p.541}), the Noether linkage sequence, and
the Gherardelli linkage theorem.  This theory is all more or less well
known, but has not always been developed exactly as here, and it is all
essential for our work in Section 3.  In Section 3, we'll prove our main
theorem, our characterization of complete intersections of codimension 2
in $\IP^n$.  Finally, we'll discuss two examples; the first shows that
the ambient projective space cannot be replaced even by a smooth
hypersurface, and the second shows that our characterization fails in
characteristic 2.

\sct2 Gherardelli linkage

\prp1 \(Peskine--Szpiro linkage theorem)\enspace
 Let $Z$ be a Gorenstein scheme, $X\subset Z$ be a proper closed
subscheme, and $Y$ the residual scheme of $X$.  If $X$ is
Cohen--Macaulay of pure codimension $0$, then so is $Y$; furthermore,
then $X$ is also the residual scheme of $Y$.
 \pf
 Let  $\cI_{X/Z}$ and  $\cI_{Y/Z}$ denote the ideals.  Then we have
  $$\cI_{Y/Z}:=\cAnn_{\cO_Z}\cI_{X/Z}
              \lisom\cHom_{\cO_Z}(\cO_X,\,\cO_Z),\eqno\Cs1.1)$$
 where the equation holds by definition and the isomorphism is given by
evaluation at 1.

It is a basic fact (see \cite{E, 21.21, \p.538}) that, on the category
of maximal (dimensional) Cohen--Macaulay $\cO_Z$-modules $\cM$, the
functor,
        $$\cD(\cM):=\cHom_{\cO_Z}(\cM,\,\cO_Z),$$
 is dualizing.  Now, $\cD$ interchanges the
two basic exact sequences,
        $$0\to\cI_{X/Z}\to\cO_Z\to\cO_X\to0\and
          0\to\cI_{Y/Z}\to\cO_Z\to\cO_Y\to0;$$
 indeed, $\cD$ carries the first sequence to the second thanks to
\Cs1.1), and so, as $\cD$ is dualizing, $\cD$ carries the second
sequence back to the first.  Thus, $\cO_Y=\cD(\cI_{X/Z})$ and
$\cI_{X/Z}= \cD(\cO_Y)$.  The latter equation implies that $X$ is the
residual scheme of $Y$.  The former equation implies that $\cO_Y$ is
maximal Cohen--Macaulay, because $\cI_{X/Z}$ is so since, at any $x\in
X$,
   $$\depth\cI_{X/Z,x}\ge\min\,(\depth\cO_{Z,x},\,1+\depth\cO_{X,x})$$
 (see \cite{E, 18.6b, \p.451}).  The proof is now complete.

\stp2   Let $P$ be a complete scheme defined over an algebraically closed
field of arbitrary characteristic.  Assume that $P$ is Gorenstein of
pure dimension at least 2, and equip $P$ with an invertible sheaf
$\cO_P(1)$, which is not necessarily ample.  For $i=1,\,2$ let $f_i\in
H^0(\cO_P(m_i))$ be a section, and $F_i:f_i=0$ its scheme of zeros.  Set
        $$\textstyle Z:=F_1\bigcap F_2,$$
 and assume that $Z$ has pure  codimension 2.

Let $X\subset Z$ be a proper closed subscheme, and assume that $X$ is
Cohen--Macaulay of pure codimension 2 in $P$.  Let $Y\subset Z$ be the
residual scheme of $X$.  By the Peskine--Szpiro linkage theorem (2.1),
also $Y$ is Cohen--Macaulay of pure codimension 2 in $P$, and $X$ is
also the residual scheme of $Y$.

\prp3 \(Noether linkage sequence)\enspace In the setup of \(2.2), the
dualizing sheaves and the ideals in $P$ are related by the following
short exact sequence:
  $$0\to\cI_{Z/P}\ox\bigomega_P(m_1+m_2)
        \to\cI_{Y/P}\ox\bigomega_P(m_1+m_2)
                \to\bigomega_X\to0.\eqno\Cs3.1)$$
 \pf First, note the following two equations:
        $$\bigomega_Z=\bigomega_P(m_1+m_2)|_Z
           \and \bigomega_X=\cI_{Y/Z}\ox\bigomega_Z.\eqno\Cs3.2)$$
 The first equation is standard, and results from basic duality theory
(see \cite{AK, Ch.~1} for example):
        $$\bigomega_Z=\cExt^2_P(\cO_Z,\,\bigomega_P)
          =\cHom_Z(\det(\cI_{Z/P}/\cI_{Z/P}^2),\,\bigomega_P|Z).$$
 The second equation in \Cs3.2) results from a series of three other
equations:
        $$\bigomega_X=\cHom(\cO_X,\,\bigomega_Z)
        =\cHom(\cO_X,\,\cO_Z)\ox\bigomega_Z=\cI_{Y/Z}\ox\bigomega_Z.$$
  These hold by elementary duality theory, by the invertiblity of
$\bigomega_Z$, and by \Cs1.1) above.

Finally, the Noether linkage sequence \Cs3.1) results from the basic
sequence,
        $$0\to\cI_{Z/P}\to\cI_{Y/P}\to\cI_{Y/Z}\to0,$$
 by tensoring it with $\bigomega_P(m_1+m_2)$ and then using the two
equations in \Cs3.2).

\rmk4 According to Enriques \cite{EC, vol.~3, p.~534}, Noether 
obtained the preceding proposition in the special case where $P$ is the
projective 3-space, Noether having stated it virtually as follows:

{\it If the curve $X$ is the partial intersection of two surfaces $F_1$
and $F_2$ of degrees $m_1$ and $m_2$, meeting further in a curve $Y$,
then the surfaces of degree $m_1+m_2-4$ passing through $Y$ cut on $X$
the complete canonical series.}

To derive this statement, take \Cs3.1), replace $\bigomega_P$ by
$\cO_P(-4)$, and extract cohomology, obtaining the following exact
sequence:
        $$H^0(\cI_{Y/P}(m_1+m_2-4))\to H^0(\bigomega_X)\to
        H^1(\cI_{Z/P}(m_1+m_2-4)).$$
 The third term vanishes as $Z$ is a complete intersection, and
Noether's statement follows.

\thm5 \(Gherardelli linkage)\enspace
 Preserve the setup of \(2.2).  Let $m_3>0$.  If an $f_3\in H^0(\cO_P(m_3))$
exists such that $Y=F_1\bigcap F_2\bigcap F_3$ where $F_3:f_3=0$, then
 $$\bigomega_X=\bigomega_P(m_1+m_2-m_3)|_X.$$
 The converse holds if, in addition,
 $$H^1(\cO_P(m_3-m_1))=0,\ H^1(\cO_P(m_3-m_2))=0,\and
 H^2(\cO_P(m_3-m_1-m_2))=0.$$
 \pf Assume an $f_3$ exists.  Then $Y=Z\bigcap F_3$.  Hence,
multiplication by $f_3$ gives a surjection
        $\mu\:\cO_Z(-m_3)\onto\cI_{Y/Z}$.
 Its kernel $\cAnn\cI_{Y/Z}(-m_3)$ is equal to $\cI_{X/Z}(-m_3)$ because
$X$ is also the residual scheme of $Y$ thanks to \Cs1).  So $\mu$
induces an isomorphism $\cO_X(-m_3)\risom\cI_{Y/Z}$.  Hence
$\bigomega_X$ has the asserted form thanks to \Cs3.2).

Conversely, assume $\bigomega_X=\bigomega_P(m_1+m_2-m_3)|_X$.  Then
twisting the Noether linkage sequence \Cs3.1) yields the following exact
sequence:
  $$0\to\cI_{Z/P}(m_3)\to\cI_{Y/P}(m_3)\to\cO_X\to0.\eqno\Cs4.1)$$
 Extracting cohomology yields the next exact sequence:
        $$H^0(\cI_{Y/P}(m_3))\to H^0(\cO_X)\to H^1(\cI_{Z/P}(m_3)).$$

Assume the additional vanishing conditions.  Then
$H^1(\cI_{Z/P}(m_3))=0$ thanks to the twisted Koszul resolution,
        $$0\to\cO_P(m_3-m_1-m_2)\to\cO_P(m_3-m_1)\textstyle
        \bigoplus\cO_P(m_3-m_2)\to\cI_{Z/P}(m_3)\to0.$$
 Hence, we may lift $1\in H^0(\cO_X)$ to an $f_3\in
H^0(\cI_{Y/P}(m_3))$.  Set $F_3:f_3=0$.

In \Cs4.1), we may replace $\cO_X$ by $\cI_{Y/Z}(m_3)$.  Hence
$\cI_{Y/Z}(m_3)$ is generated by the image of $f_3$ in
$H^0(\cI_{Y/Z}(m_3))$.  Therefore, $Y=Z\bigcap F_3$, and the proof is
complete.

\sct3 Complete intersections

\dfn1 Let $P$ be a Gorenstein scheme, $X$ a closed Cohen--Macaulay
subscheme.  We'll say that $X$ is {\it subcanonical\/} in $P$ if $P$ is
equipped with an invertible sheaf $\cO_X(1)$, and if, for some integer
$\alpha$, we have
        $$\bigomega_X=\bigomega_P(\alpha)|_X.$$
 
Assume $P$ has pure dimension at least 3, and $X$ has pure codimension
2.  We'll say that $X$ is {\it self-linked\/} in $P$ by two effective
Cartier divisors $F_1$ and $F_2$ if they meet properly in a subscheme
$Z$ containing $X$, and if $X$ is equal to the residual scheme $Y$ of
$X$ in $Z$.

\thm2 Let $P$ be a projective space of dimension $n\ge3$ over an
algebraically closed field of characteristic $0$.  Let $X\subset P$ be a
closed subscheme that is Cohen--Macaulay of pure codimension $2$.
Assume $X$ subcanonical and self-linked.  Then $X$ is a complete
intersection.

In fact, say $X$ is self-linked by hypersurfaces $F_1$ and $F_2$ of
degrees $m_1$ and $m_2$.  Then, after $F_1$ and $F_2$ are switched
if need be, $m_2$ is even, and there is a hypersurface $F_3$ of degree
$m_2/2$ such that
      $\textstyle X=F_1\bigcap F_3 \and Z=F_1\bigcap 2F_3
         \where Z:=F_1\bigcap F_2$.
 \pf
 Since $P$ is smooth and $X$ is subcanonical, $X$ is Gorenstein.  Hence,
since $X$ has pure codimension $2$, it is locally a complete
intersection in $P$ by one of Serre's results \cite{E, 21.10, p.~537}.
Hence, on $X$, the conormal sheaf $\cI_{X/P}/\cI_{X/P}^2$ is locally
free of rank 2.

By another celebrated theorem of Serre's, $H^i(\cO_P(j))=0$ for $i=1,2$
and for any $j$ since $n\ge3$.  Hence, by the Gherardelli linkage theorem
(2.5), there is a hypersurface $F_3$ such that $X=Z\bigcap F_3$.

Let $x\in X$.  For $i=1,2,3$ let $\vf_i\in\cO_{P,x}$ generate the ideal
of $F_i$.  Then $\cI_{X/P,x}$ is generated by $\vf_1$, $\vf_2$, and
$\vf_3$, but not by $\vf_1$ and $\vf_2$, since $X=Z\bigcap F_3$, but
$X\not= Z$.  Since $\cI_{X/P,\,x}$ is generated by two elements, it must
be generated either by $\vf_1$ and $\vf_3$ or by $\vf_2$ and $\vf_3$.
Hence $X$ is a Cartier divisor on $F_3$.

For $i=1,2$ set $Z_i:=F_i\bigcap F_3$.  Let $x\in X$.  Then, by the
preceding paragraph, $\cI_{X/P,\,x}$ is equal either to
$\cI_{Z_1/P,\,x}$ or to $\cI_{Z_2/P,\,x}$.  Put geometrically, $X$ is
equal, in a neighborhood of $x$ in $P$, either to $Z_1$ or to $Z_2$.

For $i=1,2,3$ say $F_i:f_i=0$.  For $i=1,2$ form the greatest common
divisor $g_i$ of $f_i$ and $f_3$, and set $G_i:g_i=0$.

First, suppose both $G_1$ and $G_2$ are nonempty, and let $x$ be a
common point.  Since $G_1$ is a component of both $F_1$ and $F_3$, their
intersection $Z_1$ is not equal to $X$ in a neighborhood of $x$.
Similarly, $Z_2$ is not equal to $X$ in a neighborhood of $x$.  This
conclusion stands in contradiction to our conclusion above, that $X$ is
equal, in a neighborhood of $x$ in $P$, either to $Z_1$ or to $Z_2$.
Therefore, not both $G_1$ and $G_2$ are nonempty; say $G_2$ is empty.

Then $Z_2$ has pure codimension 2 in $P$, and $Z_2\supseteq X$.  If
$Z_2=X$, then $X=F_2\bigcap F_3$.  So suppose not, and we'll prove that
$X=F_1\bigcap F_3$.  Form the residual scheme $X_2$ of $X$ in $Z_2$.  By
general principles, $X_2$ is a Cartier divisor on $F_3$ because $X$ and
$Z_2$ are so; moreover, $Z_2=X+X_2$.

Suppose $G_1$ is nonempty.  Set $C:=G_1\bigcap F_2$.  Then $C$ is a
hypersurface section of $F_2$.  So $C$ has a point $x$ in common with
$X_2$, which also lies on $F_2$.  Then $x\in X$, because $C\subset Z$
and $Z$ has the same support as $X$.  Since $G_1$ is a component of both
$F_1$ and $F_3$, their intersection $Z_1$ is not equal to $X$ in a
neighborhood of $x$.  Since $x$ lies on both $X_2$ and $X$, also $Z_2$
is not equal to $X$ in a neighborhood of $x$.  As before, there is a
contradiction.  Therefore, $G_1$ is empty.

Then $Z_1$ has pure codimension 2 in $P$, and $Z_1\supseteq X$.  If
$Z_1=X$, then $X=F_1\bigcap F_3$ as claimed.  So suppose not, and form
the residual scheme $X_1$ of $X$ in $Z_1$.  By general principles, $X_1$
too is a Cartier divisor on $F_3$.  After a bit of work, we'll achieve a
contradiction.

First, we'll construct a natural splitting of the natural surjection,
  $$\cI_{X/P}/\cI_{X/P}^2\onto\cI_{X/F_3}/\cI_{X/F_3}^2.\eqno\Cs2.1)$$
 To do so, in $\cI_{X/P}/\cI_{X/P}^2$ , form the image $\cL$ of
$\cI_{Z/P}$; we are going to show that $\cL$ maps isomorphically onto
$\cI_{X/F_3}/\cI_{X/F_3}^2$.  Since $\cL$ maps surjectively, and since
$\cI_{X/F_3}/\cI_{X/F_3}^2$ is invertible as $X$ is a Cartier divisor on
$F_3$, we need only show that $\cL$ is invertible.

Let $x\in X$.  Say, as above, that $\cI_{X/P,\,x}=\cI_{Z_1/P,\,x}$.  Set
$W:=F_1\bigcap 2F_3$.  Then $W\supseteq Z$; indeed,
$\cI_{F_3/P}^2\subset\cI_{Z/P}$ because $\cI_{X/Z}=\cAnn\cI_{X/Z}$ since
$X$ is self-linked.  Since also $W\supseteq Z_1$, there is a natural
commutative diagram,
 $$\CD
        0@>>>\cI_{Z_1/W}@>>>\cO_W@>>>\cO_{Z_1}@>>>0\\
        @.      @VuVV         @VvVV    @VwVV\\
        0@>>> \cI_{X/Z} @>>>\cO_Z@>>>  \cO_X  @>>>0
  \endCD$$
  Clearly, $\cI_{Z_1/W}=\cO_{P}(-F_3)|_{Z_1}$.  Moreover, $\cI_{X/Z}=
\bigomega_X\ox\bigomega_P(m_1+m_2)^{-1}$ thanks to (2.3.2) with $Y:=X$.
Thus the source of $u$ is invertible on $Z_1$, and the target is
invertible on $X$.  Now, $\cI_{X/P,\,x}=\cI_{Z_1/P,\,x}$.  Hence, $w$ is
an isomorphism at $x$; in other words, $X$ and $Z$ are the same scheme
in a neighborhood of $x$.  Also, $u$ is surjective at $x$, and its
source and target are invertible sheaves on the same scheme in a
neighborhood of $x$; hence, $u$ is an isomorphism at $x$.  So $v$ is an
isomorphism at $x$.  Hence, $\cI_{W/P,\,x}=\cI_{Z/P,\,x}$.

So, in $\cI_{X/P}/\cI_{X/P}^2$, the images of $\cI_{W/P}$ and
$\cI_{Z/P}$ are equal at $x$.  The image of $\cI_{W/P}$ is equal to
$\cO_P(-F_1)|_X$ at $x$; indeed, the latter sheaf maps naturally into
the former, and this map is surjective since $X\subset F_3$, and
injective at $x$ since its natural image is a direct summand of
$\cI_{X/P}/\cI_{X/P}^2$ at $x$, as $\cI_{X/P,\,x}=\cI_{Z_1/P,\,x}$.
  The image of $\cI_{Z/P}$ is $\cL$ by definition.  Thus $\cL$ is
invertible at $x$.  Since $x\in X$ is arbitrary, $\cL$ is invertible.
Thus $\cL\risom \cI_{X/F_3}/\cI_{X/F_3}^2$, and \Cs2.1) splits.

Let $F$ be any irreducible component of $F_3$, and equip $F$ with its
reduced structure.  Since $F$ is a hypersurface, $F$ meets $X$.  Set
$V:=X\bigcap F$.  Then $V$ is a Cartier divisor on $F$, and hence $V$ is
locally a complete intersection in $P$.  Consider the natural
commutative diagram of sheaves on $V$,
 $$\CD
        \cL|V          @>>> \bigl(\cI_{X/F_3}/\cI_{X/F_3}^2\bigr)|V\\
         @VVV                               @VVV\\
 \cI_{V/P}/\cI_{V/P}^2 @>>>         \cI_{V/F}/\cI_{V/F}^2
  \endCD$$
 The top horizontal map is an isomorphism because it is the restriction
of an isomorphism.  The right vertical map is an isomorphism because it
is surjective and its source and target are invertible.  Therefore, the
lower horizontal map splits.

Since the lower map splits, since $V$ is a Cartier divisor on $F$ and is
locally a complete intersection in $P$, since $F$ is reduced,
irreducible, and closed, and since $P$ is a projective space of
dimension $n\ge3$ over an algebraically closed field of characteristic
$0$, Braun's main theorem \cite{B, p.~26} implies that some multiple of
$V$ is numerically equivalent to a hypersurface section of $F$.

Since $F$ is a hypersurface, $F$ meets both $X_1$ and $X_2$, which are
supposedly nonempty.  For $i=1,\,2$ set $V_i:=X_i\bigcap F$.  Then $V_i$
is a Cartier divisor on $F$, and $V+V_i=F_i\bigcap F$.  Hence, some
multiple of $V_i$ too is numerically equivalent to a hypersurface
section of $F$.  Therefore, $V_1$ and $V_2$ have a common point $x$.
Then $x$ lies on both $Z_1$ and $Z_2$, so on their intersection, which
is $X$.  However, there is no neighborhood of $x$ in which either $Z_1$
or $Z_2$ is equal to $X$ because $x$ lies on both $X_1$ and $X_2$.
Thus, we've achieved the desired contradiction.  Therefore,
$X=F_1\bigcap F_3$.

Then $W=Z$ everywhere, by the reasoning above; in other words,
$Z=F_1\bigcap 2F_3$.  Finally, set $m_3:=\deg F_3$.  Then $\deg
Z=2m_1m_3$.  Now, $Z:=F_1\bigcap F_2$, so $\deg Z=m_1m_2$.  Hence
$2m_3=m_2$.  The proof is now complete.

\eg3 Most of the proof of Theorem \Cs2) works without change in the
relative case where $P$ is a smooth projectively Cohen--Macaulay variety
of pure dimension at least 3.  However, to apply Braun's theorem, we
must know that the surjection \Cs2.1) splits when $P$ is replaced by the
ambient projective space; the proof shows that \Cs2.1) itself splits,
but this splitting is insufficient.  The theorem does not hold even when
$P$ is replaced by a smooth hypersurface, as the following example
shows.

Let $P$ be a smooth quadric hypersurface in $\IP^4$.  Let $F_1$ be the
section of $P$ by a hyperplane $H_1$ that is tangent to $P$ at a point
$x$.  Then $F_1$ is a cone in $H_1$ with vertex at $x$ and with base a
smooth (plane) conic $C$.  Fix $y\in C$.  Then $y$ determines a
generator $X$ of the cone $F_1$.  Let $H_2$ be a hyperplane in $\IP^4$
that cuts $H_1$ in the plane spanned by $x$ and by the tangent line to
$C$ at $y$.  Then $X$ is a line, so subcanonical in $P$.  Moreover, $X$
is self-linked in $P$ by $F_1$ and $F_2$ with $F_2:=H_2\bigcap P$.
However, $X$ is not the complete intersection of two hypersurface
sections of $P$ since any such complete intersection has even degree in
$\IP^4$.

\eg4 Theorem \Cs2) is not valid in positive characteristic without some
further restriction on $X$.  Indeed, we are going to see that, in
characteristic 2, there exists an example of an irreducible, but
nonreduced, Cohen--Macaulay space curve $X$, which is subcanonical and
self-linked, yet is not a complete intersection.

Ferrand \cite{F, p.~345} explained how to put a subcanonical double
structure on a line (indeed, on any complete curve that is locally a
complete intersection) in $\IP^3$ in any characteristic; moreover, the
double curve can have arbitrarily negative arithmetic genus.  Now,
Migliore \cite{M, p.~185} proved that, in characteristic 2, a double
line $X$ is self-linked if its arithmetic genus is $-2$ or less.  Such
an $X$ is not a complete intersection, because every complete
intersection $Z$ has nonnegative arithmetic genus by (2.3.2).

 \references
 \serial{AUF}{Ann. Univ. Ferrara}
 \serial{annali}{Annali di Mat.}
 \serial{ArchMath}{Arch. Math.}
 \serial{cras}{C. R. Acad. Sc. Paris}
 \serial{Duke}{Duke Math. J.}
 \serial{inv}{Inventiones math.}
 \serial{rsmupt}{Rend. Sem. Univ. Pol. Torino}
 \serial{RendiMilano}{Rend. Sem. Mat. Fis. Milano}
 \serial{RendiAccIt}{Rend. Acc. d'Italia}
 \serial{tams}{Trans. Amer. Math. Soc.}

AK
 A. Altman and S. Kleiman,
 ``Introduction to Grothendieck duality theory,''
 Lecture Notes in Math. {\bf 146}, Springer-Verlag, 1970.

BP
 B. Basili and C. Peskine,
 D\'ecomposition du fibr\'e normal des surfaces lisses de  $\IP_4$ et
structures doubles sur les solides de $\IP_5$, \Duke 69 1993 87--95

BE
 V. Beorchia and Ph. Ellia,
 Normal bundle and complete intersections,
 \rsmupt 48 1990 553--62

B
 R. Braun,
 On the normal bundle of Cartier divisors on projective varieties,
 \ArchMath 59 1992 403--11

E
 D. Eisenbud,
 ``Commutative Algebra,''
 GTM {\bf 150}, Springer-Verlag 1994.

EC
 F. Enriques and O. Chisini,
 ``Lezioni sulla teoria geometrica delle equazioni e delle funzioni
algebriche,'' Zanichelli, Bologna, 1915.

Fa
 G. Faltings,
 Ein Kriterium for vollst\"andige Durchshnitte,
 \inv 62 1980 393--401

F
 D. Ferrand,
 Courbes gauches et fibr\'es de rang $2$,
 \cras 281 1975 345--7

FL81
 M. Fiorentini and A. T. Lascu,
 Una formula di geometria numerativa,
 \AUF 27 1981 201--27

FL87
 M. Fiorentini and A. T. Lascu,
 Projective embeddings and linkage,
 \RendiMilano LXVII 1987 161--82

Gh
 G. Gherardelli,
 Sulle curve sghembe algebriche intersezioni complete di tre superficie,
 \RendiAccIt IV 1943 460--62

M
 J. Migliore,
 On linking double lines,
 \tams 294 1986 177--85

PS
 C. Peskine and L. Szpiro,
 Liaison des vari\'et\'es alg\'ebriques. I,
 \inv 26 1974 271--302

R79
 P. Rao,
 Liaison among curves in $\IP^3$,
 \inv 50 1979 205--17

R82
 P. Rao,
 On self-linked curves,
 \Duke 49 1982 251--73

\endreferences
\bye